\documentclass[11pt]{amsart}

\usepackage{amsmath}
\usepackage{amsfonts}
\usepackage{amssymb}
\usepackage{amscd}
\usepackage{amsthm}
\usepackage{framed}
\usepackage{fullpage}
\usepackage{graphicx}
\usepackage{latexsym}
\usepackage[numbers]{natbib}  

\usepackage{hyperref}
\hypersetup{colorlinks,linkcolor={red},citecolor={blue},urlcolor={blue}}

\newtheorem{theorem}{Theorem}[section]
\newtheorem{lemma}[theorem]{Lemma}
\newtheorem{proposition}[theorem]{Proposition}
\newtheorem{corollary}[theorem]{Corollary}

\theoremstyle{definition}

\newtheorem{problem}[theorem]{Problem}
\newtheorem{conjecture}[theorem]{Conjecture}

\numberwithin{equation}{section}

\newcommand{\Z}{\mathbb Z}
\newcommand{\Q}{\mathbb Q}
\newcommand{\R}{\mathbb R}
\newcommand{\C}{\mathbb C}
\newcommand{\F}{\mathbb F}
\newcommand{\OO}{\mathrm O}
\newcommand{\SO}{\mathrm{SO}}
\newcommand{\PO}{\mathrm{PO}}
\newcommand{\D}{\mathcal D}
\newcommand{\HH}{\mathcal H}
\newcommand{\vol}{\operatorname{vol}}
\newcommand{\HM}{\mathrm{HM}}
\newcommand{\rk}{\operatorname{rk}}
\newcommand{\im}{\operatorname{im}}
\newcommand{\divisor}{\operatorname{div}}
\newcommand{\Fix}{\operatorname{Fix}}
\newcommand{\fpr}{\operatorname{fpr}}
\newcommand{\Sub}{\operatorname{Sub}}

\newcommand{\Pic}{\operatorname{Pic}}
\newcommand{\Cl}{\operatorname{Cl}}

\newcommand{\cont}{\operatorname{cont}}
\newcommand{\eps}{\varepsilon}

\begin{document}

\title[Finiteness of strongly reflective modular forms]{Finiteness of strongly reflective modular forms} 

\author{Jialin Li}

\address{School of Mathematics and Statistics, Wuhan University, Wuhan 430072, Hubei, China}

\email{jlli.math@whu.edu.cn}

\author{Haowu Wang}

\address{School of Mathematics and Statistics, Wuhan University, Wuhan 430072, Hubei, China}

\email{haowu.wangmath@whu.edu.cn}


\date{\today}


\begin{abstract}
Let $\Gamma\subset\mathrm{O}^+(n,2)$ be a rational arithmetic group with $n\geq 3$. A non-constant modular form for $\Gamma$ is called strongly reflective if its zero divisor is a reduced union of mirrors of reflections in $\Gamma$. We show that, up to conjugation and multiplication by a nonzero scalar, there are only finitely many strongly reflective modular forms. As a consequence, we establish the finiteness of holomorphic automorphic products of singular weight. These results answer a question posed by Borcherds in 1995 and settle a conjecture of Gritsenko–Nikulin from 1998. As an application, up to conjugation, we prove the finiteness of $\Gamma$ whose corresponding algebra of modular forms is freely generated. The key ingredient is a uniform upper bound on the covolume of $\Gamma$, obtained by adapting the higher-rank lattice convergence theorem of Abért et al. to the study of the degrees of reflective branch divisors.  
\end{abstract}

\maketitle

\section{Introduction}
Let $\Gamma$ be a rational arithmetic subgroup of $\mathrm{O}(n,2)$, where $n\geq 3$. Realize $\Gamma$ as a finite-index subgroup of the integral orthogonal group $\mathrm{O}(L)$, where $L$ is an even integral lattice of signature $(n,2)$. Put $V=L\otimes\Q$, $V_\R=V\otimes_\Q\R$, and $V_\C=V\otimes_\Q\C$. The associated symmetric domain $\D_L$ of type IV is a connected component of the set
\[
 \{[Z]\in\mathbb P(V_\C):\; (Z,Z)=0,\ (Z,\overline Z)<0\}.
\]
Denote by $\OO^+(L)$ the subgroup of lattice isometries that preserve $\D_L$ and assume $\Gamma \subset \mathrm{O}^+(L)$. Its action on $\D_L$ has kernel $\{\pm1\}$, so we often consider the projective group
\[
\PO^+(L)=\OO^+(L)/\{\pm1\},
\]
which acts faithfully on $\D_L$. The affine cone over $\D_L$ is 
\[
\D_L^\bullet=\{Z\in V_\C\setminus\{0\}:\; [Z]\in\D_L\}. 
\]
For an integer $k$ and a character $\chi:\Gamma\to\C^\times$, a \emph{holomorphic modular form} of weight $k$ for $\Gamma$ and character $\chi$ is a holomorphic function $F:\D_L^\bullet\to\C$ that satisfies
\[
\begin{aligned}
F(tZ)&=t^{-k}F(Z),&& \text{for any $t\in\C^\times$},\\
F(\gamma Z)&=\chi(\gamma)F(Z),&& \text{for any $\gamma\in \Gamma$}. 
\end{aligned}
\]
This definition extends naturally to half-integral weights, in which case $\chi$ becomes a multiplier system. The weight of any non-constant $F$ satisfies
the lower bound $k\geq n/2-1$; this minimal possible value is called the \textit{singular weight}. 

In 1995 and 1998, Borcherds \cite{Bor95, Bor98} developed the singular theta lift, which constructs modular forms on $\mathrm{O}(n,2)$ with infinite product expansions. Such modular forms are usually candidates for denominator functions of generalized Kac--Moody algebras, which played a central role in Borcherds' celebrated proof of the moonshine conjecture \cite{Bor92}. 

The singular theta lift is multiplicative. Its input $f$ is a vector-valued modular form of weight $1-n/2$ for the Weil representation of $\mathrm{Mp}_2(\Z)$ attached to $L$, holomorphic on the upper half-plane and meromorphic at the cusps, with integral Fourier coefficients of nonpositive index. The output $\Psi(f)$ is a meromorphic modular form of weight $c_0(0)/2$ for a finite-index subgroup of $\mathrm{O}^+(L)$ with a finite-order multiplier system, where $c_0(0)$ is the constant coefficient of the zero component of $f$. These outputs are called \emph{Borcherds products} or \emph{automorphic products}. They have infinite product expansions at each zero-dimensional cusp, and their zeros and poles lie on the \emph{rational quadratic divisors} $\D_v$ with multiplicities determined by Fourier coefficients of negative index of $f$, where
\[
\D_v=\{[Z]\in\D_L: \; (Z,v)=0\},\quad v\in V,\quad(v,v)>0.
\]
The divisor $\D_v$ is the mirror of the orthogonal reflection 
\[
\sigma_v:\; V_\R \to V_\R, \quad x \mapsto x-\frac{2(x,v)}{(v,v)}v.
\]
We write $r_v\in \PO^+(n,2)$ for its projective action. When $\Psi(f)$ is holomorphic, it is a holomorphic modular form in the sense defined above.

Holomorphic automorphic products of singular weight are of particular interest: their Fourier coefficients are supported only on isotropic vectors, and the corresponding generalized Kac--Moody superalgebras may admit natural constructions with known symmetries and have deep connections with vertex algebras and affine Lie superalgebras \cite{Bor90, Bor92, Sch00, SWW23, DSW23, SWW26}. In 1995, Borcherds posed the following problem:

\begin{problem}[{\cite[Problem 3, Section 17]{Bor95}}]\label{problem:Bor95}
Are there a finite or infinite number of automorphic forms of singular weight that can be written as modular products? 
\end{problem} 

The first substantial step toward this problem was taken by Williams and the second named author \cite{WW23}. For a singular automorphic product $F$ on $L$, it was proved that $F$ has only simple zeros and that, whenever $F$ vanishes along $\D_v$, one has $F(\sigma_v(Z))=-F(Z)$. As a consequence, $F$ can be regarded as a strongly reflective modular form.

The notion of reflective modular forms was introduced by Borcherds \cite{Bor98} and Gritsenko--Nikulin \cite{GN98a, GN98}. A non-constant holomorphic modular form for $\Gamma\subset \mathrm{O}^+(L)$ is called \emph{reflective} if its zero divisor is a linear combination of the $\D_v$ with $\sigma_v\in \Gamma$; it is further called \emph{strongly reflective} if the multiplicity of each $\D_v$ is at most one. By Bruinier's converse theorem \cite{Bru02, Bru14}, reflective modular forms can often be realized as automorphic products. Motivated by their study of automorphic corrections of Lorentzian Kac–Moody algebras, Gritsenko and Nikulin \cite{GN98, GN98a} conjectured in 1998 the finiteness of strongly reflective modular forms. They focused on the simple-zero case, since some zero mirrors correspond to real roots of generalized Kac–Moody Lie algebras, and such roots have multiplicity one. We reformulate their conjecture as follows: 

\begin{conjecture}[{\cite[Conjecture 2.5.5]{GN98a}}]\label{conjecture:GN98}
There are only finitely many lattices that carry strongly reflective modular forms, and only finitely many such modular forms, up to certain equivalences. 
\end{conjecture}

We now explain why every singular automorphic product $F$ is strongly reflective. As shown in \cite[Lemmas 3.1--3.2]{WW23}, if one adjoins the zero reflections of $F$ to its original modular group and embeds the extended modular group into the integral orthogonal group of $\mathbb{L}$, then $F$ becomes a strongly reflective modular form on $\mathbb{L}$, where $\mathbb{L}$ is typically constructed by rescaling a suitable sublattice of $V$. Consequently, Conjecture \ref{conjecture:GN98} establishes finiteness, answering Problem \ref{problem:Bor95}. 

Reflective modular forms have found wide applications in infinite-dimensional Lie algebras, algebraic geometry, hyperbolic reflection groups, and free algebras of modular forms \cite{Bor98, GN98a, GHS07b, GH14, Ma18, Wan21a}. Their classification has been an active area of research in the past 30 years \cite{GN98, GN02, Sch06, Sch17, Ma17, Ma18, Dit19, Wan25b, SWW23, DSW23, SWW26}. However, most of these results focus on the classification of the underlying lattices, and the finiteness of reflective modular forms remains open in full generality since 1998. The main difficulty lies in bounding the modular groups $\Gamma$ of reflective modular forms when $\Gamma$ varies over finite-index subgroups of a fixed integral orthogonal group $\mathrm{O}^+(L)$. 

We overcome this difficulty by establishing a uniform upper bound on the covolume of $\Gamma$ and then applying Wang's finiteness theorem for lattices of bounded covolume \cite[Theorem 8.1]{Wan83}. We derive this bound by expressing the degrees of reflective branch divisors in terms of fixed-point ratios associated with invariant random subgroups and invoking the higher-rank lattice convergence theorem of Ab\'ert et al. \cite{ABBGNRS17, AGN17}. 

The following is our main theorem, which completely solves Problem \ref{problem:Bor95} and Conjecture \ref{conjecture:GN98}.

\begin{theorem}\label{Main theorem}
For any fixed positive integer $N$, there are only finitely many reflective modular forms whose zero mirrors all have multiplicity at most $N$, up to the equivalence
\[
F_1\sim F_2 \Longleftrightarrow \text{there exist $c\in\C^\times$ and $\gamma\in\mathrm{O}(n,2)$ such that $F_2=c\cdot(F_1|_k\gamma)$}.
\]
In particular, up to the same equivalence, both strongly reflective modular forms and holomorphic automorphic products of singular weight are finite in number. 
\end{theorem}

The above finiteness result holds in full generality: we allow all weights $k\in\frac{1}{2}\Z$ and all multiplier systems, consider all finite-index subgroups of $\mathrm{O}^+(L)$, and let $L$ range over all even integral lattices of signature $(n,2)$ for any integer $n\geq 3$. 

We present a further application of the main theorem. Let $M_*(\Gamma)$ be the algebra of modular forms of integral weights with trivial character for $\Gamma\subset \mathrm{O}(n,2)$, where $n\geq 3$. This algebra is an infinite counterpart of the polynomial invariants of a finite linear group. Note that the freeness of $M_*(\Gamma)$ implies that $\Gamma$ is generated by reflections and that $n\leq 10$ (see \cite{SV17} and \cite[Theorem 1.3]{Wan25b}). In analogy with the seminal Shephard–Todd–Chevalley theorem, the second named author \cite{Wan21a} proved that $M_*(\Gamma)$ is freely generated if and only if the Jacobian of potential free generators is a strongly reflective modular form for $\Gamma$ that vanishes precisely along all mirrors of reflections in $\Gamma$. Combining this criterion and Theorem \ref{Main theorem}, we prove the finiteness of $\Gamma$ for which $M_*(\Gamma)$ is free. 

\begin{corollary}
Up to conjugation, there are only finitely many rational arithmetic groups $\Gamma\subset \mathrm{O}^+(n,2)$ whose corresponding algebra $M_*(\Gamma)$ of modular forms is freely generated over $\C$.     
\end{corollary}

We now outline the proof of Theorem \ref{Main theorem}. Fix an even integral lattice $L$ of signature $(n,2)$ and put $\Lambda=\PO^+(L)$. Let $\Gamma\subset \Lambda$ be a finite-index subgroup and denote by $\widetilde\Gamma$ its full inverse image in $\OO^+(L)$. Assume that $F$ is a reflective modular form of weight $k$ for $\widetilde\Gamma$, with the multiplicities of its zero mirrors bounded by $N$. On the Baily--Borel compactification \cite{BB66}, let $\lambda_\Gamma$ be the rational automorphic class and let $B_F$ record the images of the zero mirrors of $F$ with their multiplicities on $\D_L$. Then the normalized degree $T(\Gamma)$ of the reduced branch divisor $B_\Gamma$ for the quotient map $\D_L \to \Gamma \backslash \D_L$ satisfies (see Lemma \ref{lem:divisor})
\begin{equation}\label{eq:step 1}
n-2\leq 2k=\frac{B_F\cdot\lambda_\Gamma^{n-1}}{\lambda_\Gamma^n}\leq N\cdot T(\Gamma), \qquad T(\Gamma):=\frac{B_\Gamma\cdot\lambda_\Gamma^{n-1}}{\lambda_\Gamma^n}. 
\end{equation}
Moreover, we derive the expression (see Proposition \ref{prop:branch-cover})
\begin{equation}\label{eq:step 2}
T(\Gamma)=\sum_{j=1}^s c_j\cdot \fpr_\Gamma(r_j), \qquad c_j:=\frac{B_j\cdot\lambda_\Lambda^{n-1}}{\lambda_\Lambda^n}>0, \quad \fpr_\Gamma(r_j):=\frac{|\Fix_{\Lambda/\Gamma}(r_j)|}{[\Lambda: \Gamma]},
\end{equation}
where $\Fix_{\Lambda/\Gamma}(r_j):=\{ g\Gamma \in \Lambda/\Gamma:\, r_jg\Gamma=g\Gamma \}$, and $r_1,...,r_s\in \Lambda$ are projective reflections giving rise to $\Lambda$-orbits for all mirrors of reflections in $\Lambda$, and each $B_j$ is the corresponding reduced divisor. 

By the higher-rank convergence theorems \cite{ABBGNRS17,AGN17}, for finite-index subgroups $\Gamma_i\subset \Lambda$, we have 
\begin{equation}\label{eq:step 3}
[\Lambda:\Gamma_i]\longrightarrow\infty
\quad\Longrightarrow\quad \fpr_{\Gamma_i}(r_j) \longrightarrow 0, \; \forall\, 1\leq j\leq s \quad\Longrightarrow \quad T(\Gamma_i)\longrightarrow 0.
\end{equation}
Here, we have to adapt the convergence results of \cite{ABBGNRS17,AGN17} (see Proposition \ref{prop:fpr}), since projective reflections do not belong to the identity component of $\PO^+(n,2)$ when $n$ is even. 

The lower bound $T(\Gamma)\geq (n-2)/N$, together with the above limit, yields a uniform upper bound on the indices $[\Lambda:\Gamma]$. Hence, up to conjugation, only finitely many such subgroups $\Gamma$ occur, by the finite generation of $\Lambda$ \cite{BH62} or Wang's finiteness theorem for lattices of bounded covolume \cite{Wan83}. For each fixed $\Gamma$, the finiteness of the mirror orbits, together with the multiplicity bound $N$, shows that only finitely many modular forms $F$ occur up to multiplication by a nonzero scalar. 

To complete the proof, it remains to establish the finiteness, up to isometry, of primitive integral lattices $L$ admitting a reflective modular form with vanishing orders bounded by $N$ for a finite-index subgroup $\Gamma\subset\mathrm{O}^+(L)$. Ma \cite[Corollary 1.9]{Ma18} proved this for $n\geq 4$. We adapt Ma's volume estimates to treat the remaining case $n=3$. For every $n\geq 3$, we establish (see Theorem \ref{thm:bounds})
\begin{equation}\label{eq:step 4}
\frac{n-2}{N}\leq T(\Gamma)\leq T(\Lambda)\leq C_n \cdot D(L)^{-1/4},
\end{equation}
where $D(L)$ is the exponent of the discriminant group $L^\vee/L$, and $C_n>0$ is an effective constant depending only on $n$. Thus, $D(L)$ is bounded in terms of $n$ and $N$, so that, for each fixed $n$, there are only finitely many isometry classes of $L$. Finally, the second named author \cite[Theorem 1.2]{Wan25b} has shown that reflective modular forms on $\mathrm{O}(n,2)$ exist if and only if $n\leq 20$ or $n=22,26$. This concludes the proof of Theorem~\ref{Main theorem}.

\section{Fixed-point ratios and branch degrees}\label{sec:main step}
In this section, we recall the higher-rank lattice convergence theorem and introduce branch degrees. We establish \eqref{eq:step 1}--\eqref{eq:step 3}, which constitute the main steps in the proof of Theorem \ref{Main theorem}.  

We fix the notation as at the beginning of the introduction. Let $\OO^+(V_\R)$ be the real orthogonal group preserving the chosen domain component. The projective group
\[
 G:=\PO^+(V_\R)=\OO^+(V_\R)/\{\pm1\}
\]
acts effectively on the symmetric domain $\D_L$. An \emph{effective rational orthogonal arithmetic group} is a subgroup $\Gamma$ of
\[
\PO^+(V):=\im\bigl(\OO(V)\cap\OO^+(V_\R)\longrightarrow G\bigr)
\]
whose intersection with $\PO^+(L)$ has finite index in both groups. Here $\OO(V)$ is the group of rational isometries of $V:=L\otimes\Q$.

For $Z=X+iY\in\D^\bullet_L$, we have $(X,X)=(Y,Y)<0$ and $(X,Y)=0$. Thus, $\D_L$ parametrizes the oriented negative definite real two-planes. The identity component $G^\circ$ is the connected component of $1$ in the real Lie group $G$. It is the image of $\SO_0(V_\R)$, the identity component of the real special orthogonal group.

For subgroups $\Delta$, $\Delta'$ of a common group and an element $a$ of that group, write
\[
C_\Delta(\Delta')=\{g\in\Delta:\, gh=hg, \, \text{for every $h\in \Delta'$}\},\qquad
C_\Delta(a)=\{g\in\Delta:\, ga=ag\}. 
\]
These are the centralizers of $\Delta'$ and $a$ in $\Delta$. A lattice in a real Lie group is a discrete subgroup whose quotient has finite Haar volume. The effective rational orthogonal arithmetic groups are lattices in $G$ by Borel--Harish-Chandra \cite{BH62}. 

\begin{lemma}\label{lem:group-background}
The identity component $G^\circ$ is center-free and simple of real rank two. The group $G$ is connected when $n$ is odd and has two components when $n$ is even. Moreover,
\[
 C_G(G^\circ)=\{1\}.
\]
Every finite-index subgroup of a lattice in $G^\circ$ is Zariski dense in the adjoint real algebraic group $\mathbf{H}:=\mathbf{SO}(V_\R)_{\mathrm{ad}}\simeq \mathbf{SO}(V_\R)/Z(\mathbf{SO}(V_\R))$ with $G^\circ\simeq \mathbf{H}(\R)^\circ$. 
\end{lemma}

\begin{proof}
Simplicity, real rank, and Borel density are given, respectively, by \cite[Example A2.3 (2), Proposition 8.1.8, and Corollary 4.5.6]{Mor15}. A finite-index subgroup is again a lattice.

The Cartan decomposition gives $\OO^+(n,2)$ the same components as its maximal compact subgroup $\OO(n)\times\SO(2)$. Multiplication by $-1$ interchanges them exactly when $n$ is odd. If the image of an orthogonal lift $A$ centralizes $G^\circ$, the commutator of $A$ with an element of $\SO_0(n,2)$ lies in $\{\pm1\}$. It equals $1$ by connectedness. The complexified standard representation is irreducible, so $A$ is scalar, and orthogonality gives $A=\pm1$. This proves $C_G(G^\circ)=\{1\}$ and the center-freeness of $G^\circ$.
\end{proof}

By Lemma~\ref{lem:group-background} and \cite[Theorems 1.6.1 and 1.7.1]{BHV08}, both $G^\circ$ and its lattices satisfy the Kazhdan property~(T), meaning that every unitary representation with almost invariant unit vectors has a nonzero invariant vector. Such lattices have finite abelianization \cite[Corollary 1.3.6]{BHV08}. 

For the modular forms considered in Theorem~\ref{Main theorem}, $F^4$ has even integral weight $4k$ and an ordinary character of $\Gamma$. Squaring removes the square-root ambiguity, and even weight makes the action of $-1$ trivial. The character has finite image on $\Gamma\cap G^\circ$, hence on $\Gamma$ by finite index.

Let $H$ be a connected center-free semisimple Lie group without compact factors. We use $\vol$ for ordinary Haar or Riemannian volume, with an invariant measure fixed on each homogeneous space. The space $\Sub(H)$ of closed subgroups is compact and metrizable in its \emph{Chabauty topology}. Convergence $J_i\to J$ means that every limit of a convergent sequence $x_{i_j}\in J_{i_j}$ lies in $J$ and every $x\in J$ is a limit of elements $x_i\in J_i$. For a lattice $\Delta<H$, the map
\[
 H/\Delta\longrightarrow\Sub(H),\qquad
 h\Delta\longmapsto h\Delta h^{-1}
\]
pushes the normalized Haar measure forward to the \emph{invariant random subgroup} $\mu_\Delta$. Its invariance means invariance under conjugation by $H$. More precisely, for any Borel set $\mathcal{E}\subset \Sub(H)$, 
\[
\mu_\Delta(\mathcal{E})=\frac{\vol\{h\Delta :\, h\Delta h^{-1} \in \mathcal{E}\}}{\vol(H/\Delta)}, \qquad \text{and} \quad \mu_\Delta(\mathcal{E})=\mu_\Delta(a^{-1}\mathcal{E}a), \; \text{for any $a\in H$}. 
\]

We now recall the weak-convergence consequence of \cite[Theorem 4.4]{ABBGNRS17}. Here, weak convergence $\mu_{\Delta_i}\Rightarrow\delta_{\{1\}}$ means that $\int\varphi\,d\mu_{\Delta_i}\to\varphi(\{1\})$ for every continuous function $\varphi$ on $\Sub(H)$, and $\delta_{\{1\}}$ is the point mass at the trivial subgroup. 

\begin{theorem}\label{thm:weak-convergence}
Assume that $H$ has real rank at least two and satisfies the Kazhdan property~(T). If $\Delta_i<H$ are irreducible lattices and $\vol(H/\Delta_i)\to\infty$, then $\mu_{\Delta_i}\Rightarrow\delta_{\{1\}}$.
\end{theorem}

The theorem applies to $H=G^\circ$, and simplicity ensures that every lattice in $G^\circ$ is irreducible. In the remainder of this paper, we focus on this case.

For the finite-cover calculation below, let $\Lambda$ be a lattice in $G^\circ$ and let $\Gamma\leq \Lambda$ be a finite-index subgroup. The finite set $\Lambda/\Gamma=\{g\Gamma:\, g\in\Lambda\}$ carries the left multiplication action of $\Lambda$. Define
\[
\Fix_{\Lambda/\Gamma}(a):=\{g\Gamma:\, ag\Gamma=g\Gamma\},
\qquad
\fpr_\Gamma(a):=\frac{|\Fix_{\Lambda/\Gamma}(a)|}{[\Lambda:\Gamma]}.
\]
This is the probability that a uniformly chosen coset is fixed by $a$.

For finite-index subgroups of a fixed lattice, Ab\'ert--Gelander--Nikolov \cite{AGN17} prove a discrete version of the above weak convergence. We express it in terms of fixed-point ratios. 

\begin{lemma}\label{lem:connected-fpr}
Let $\Lambda<G^\circ$ be a lattice. If the finite-index subgroups $\Gamma_i\leq\Lambda$ satisfy $[\Lambda:\Gamma_i]\to\infty$, then $\fpr_{\Gamma_i}(a)\to0$, for every fixed $a\in\Lambda\setminus\{1\}$.
\end{lemma}

\begin{proof}
By Lemma~\ref{lem:group-background}, $G^\circ$ is a higher-rank simple Lie group. According to \cite[Theorems 10, 13 and Lemma 14]{AGN17}, we have
\[
 \nu_i:=\frac1{[\Lambda:\Gamma_i]}
       \sum_{g\Gamma_i\in\Lambda/\Gamma_i}\delta_{g\Gamma_i g^{-1}}
 \ \Longrightarrow\ \delta_{\{1\}}
 \quad\text{on }\Sub(\Lambda).
\]
The sum assigns an equal mass to each conjugate of $\Gamma_i$. For $a\ne1$, the function $J\mapsto\mathbf1_{\{a\in J\}}$ on $\Sub(\Lambda)$ is continuous and vanishes at $\{1\}$, since $\Lambda$ is discrete. We then derive
\[
 \fpr_{\Gamma_i}(a)
 =\nu_i(\{J\leq\Lambda:a\in J\})\longrightarrow0.\qedhere
\]
\end{proof}

We extend Lemma~\ref{lem:connected-fpr} to the full projective group, as required for the branch-degree estimates. When $n$ is even, both lifts $\pm\sigma_v$ of a projective reflection $r_v$ have determinant $-1$, so $r_v\notin G^\circ$.

\begin{proposition}\label{prop:fpr}
Fix an arithmetic lattice $\Lambda<G:=\PO^+(n,2)$, where $n\geq3$. If the finite-index subgroups $\Gamma_i\leq\Lambda$ satisfy $[\Lambda:\Gamma_i]\to\infty$, then
\begin{equation}\label{eq:fpr}
 \fpr_{\Gamma_i}(a)\longrightarrow 0,
 \qquad \text{for every fixed $\, a\in\Lambda\setminus\{1\}$}.
\end{equation}
\end{proposition}

\begin{proof}
Apply Lemma~\ref{lem:connected-fpr} to the identity component as follows. Put $\Lambda^\circ=\Lambda\cap G^\circ$ and $\Gamma_i^\circ=\Gamma_i\cap G^\circ$. The normal subgroup $\Lambda^\circ$ has index at most two in $\Lambda$, so its action on $\Lambda/\Gamma_i$ has at most two orbits. Note that $\mathrm{Stab}_{\Lambda^\circ}(g\Gamma_i)=g\Gamma_i^\circ g^{-1}$. Thus the orbit of $g\Gamma_i$ is in bijection with 
\[
 \Lambda^\circ/(g\Gamma_i^\circ g^{-1}).
\]
The subgroup-index formula gives
\begin{equation}\label{eq:orbit-indices}
 [\Lambda^\circ:g\Gamma_i^\circ g^{-1}]
 =[\Lambda^\circ:\Gamma_i^\circ]
 \geq\frac12[\Lambda:\Gamma_i]\longrightarrow\infty.
\end{equation}
Since $[\Lambda:\Lambda^\circ]\leq 2$, the representative $g$ of the $\Lambda^\circ$-orbit can be chosen from a fixed set of at most two representatives for $\Lambda/\Lambda^\circ$. Apply Lemma~\ref{lem:connected-fpr} to each sequence $g\Gamma_i^\circ g^{-1}<\Lambda^\circ$. The fixed-point ratio on $\Lambda/\Gamma_i$ is a weighted average of the ratios on its $\Lambda^\circ$-orbits, so
\begin{equation}\label{eq:fpr-connected-case}
 \fpr_{\Gamma_i}(a)\longrightarrow0
 \qquad(a\in\Lambda^\circ\setminus\{1\}).
\end{equation}

Now let $a\in\Lambda\setminus\Lambda^\circ$. If its $\Lambda^\circ$-conjugacy class were finite, $C_{\Lambda^\circ}(a)$ would have finite index. Borel density would then force $a$ to centralize $G^\circ$, contradicting $C_G(G^\circ)=\{1\}$ in Lemma~\ref{lem:group-background}.

Fix $m\geq1$ and choose distinct $\Lambda^\circ$-conjugates $a_1,\ldots,a_m$ of $a$. They have the same fixed-point ratio. More precisely, every $A_{i,j}:=\Fix_{\Lambda/\Gamma_i}(a_j)$ satisfies $|A_{i,j}|=[\Lambda:\Gamma_i]\cdot \fpr_{\Gamma_i}(a)$. For any $j\neq l$, 
\[
a_j^{-1}a_l\in\Lambda^\circ\setminus\{1\} \qquad \text{and} \qquad
A_{i,j}\cap A_{i,l}
\subseteq\Fix_{\Lambda/\Gamma_i}(a_j^{-1}a_l).
\]
From the fact $\bigcup_{j=1}^m A_{i,j}\subset \Lambda/ \Gamma_i$ and the elementary inequality
\[ 
\left|\bigcup_{j=1}^m A_{i,j}\right|
\ge
\sum_{j=1}^m|A_{i,j}|
-\sum_{j<l}|A_{i,j}\cap A_{i,l}|,
\]
we derive
\[
 m\,\fpr_{\Gamma_i}(a)
 \leq1+\sum_{j<l}\fpr_{\Gamma_i}(a_j^{-1}a_l).
\]
By \eqref{eq:fpr-connected-case}, the sum tends to zero for every fixed $m$, which further leads to
\[
0\leq  \limsup_i\fpr_{\Gamma_i}(a)\leq\inf_{m\geq1}\frac1m=0.\qedhere
\]
\end{proof}

We now apply fixed-point ratios to branch divisors on arithmetic quotients. For any $v\in V$ of positive norm, the reflection in the mirror $\D_v$ is
\begin{equation}\label{eq:reflection}
 \sigma_v(x)=x-\frac{2(x,v)}{(v,v)}v, \qquad x\in V_\C. 
\end{equation}
The mirror has complex dimension $n-1$. A mirror $\mathcal D_v$ is called a branch mirror for $\Gamma$ if its projective reflection $[\sigma_v]$, which fixes $\D_v$ pointwise, belongs to $\Gamma$. For simplicity, we use $\mathcal H=\D_v$ and $r=r_v=[\sigma_v]$ for a mirror and its projective reflection, respectively. For an orthogonal group with image $\Gamma$, the condition $r_v\in\Gamma$ means that the group contains $\sigma_v$ or $-\sigma_v$.

For an effective rational orthogonal arithmetic group $\Gamma$, put $X_\Gamma=\Gamma\backslash\D_L$ and write $Y_\Gamma=X_\Gamma^{\mathrm{BB}}$ for its Baily--Borel compactification \cite{BB66}. This is a normal projective variety with an ample rational automorphic class
\[
\lambda_\Gamma=\frac1a[\mathcal L_{\Gamma,a}]\in\Pic(Y_\Gamma)\otimes\Q,
\qquad
M_{aj}(\Gamma)=H^0(Y_\Gamma,\mathcal L_{\Gamma,a}^{\otimes j})
\]
for a sufficiently divisible even weight $a$. Here, $M_{aj}(\Gamma)$ is the space of weight-$aj$ forms with trivial character, and $\Pic$ is the Picard group. Write $c_1$ for the first Chern class. For a Weil divisor $E$, set
\[
 E\cdot\lambda_\Gamma^{n-1}
 =\frac{E\cdot c_1(\mathcal L_{\Gamma,a})^{n-1}}{a^{n-1}},
 \qquad
 \lambda_\Gamma^n=\frac{c_1(\mathcal L_{\Gamma,a})^n}{a^n}.
\]
Top intersections denote numerical degrees. Ampleness gives
\[
 \lambda_\Gamma^n>0,\qquad
 E\cdot\lambda_\Gamma^{n-1}>0\quad(E\ne0\text{ effective}).
\]
The Baily--Borel boundary components have complex dimension at most one, so
\[
 \dim(Y_\Gamma\setminus X_\Gamma)\leq1,\qquad
 \operatorname{codim}_{Y_\Gamma}(Y_\Gamma\setminus X_\Gamma)\geq2
 \quad(n\geq3).
\]
Consequently, every prime divisor meets $X_\Gamma$. Normality and projectivity also give
\[
 f\in\C(Y_\Gamma)^\times,\quad\divisor(f)=0
 \quad\Longrightarrow\quad f\in\C^\times.
\]
Let $p_\Gamma:\D_L\to X_\Gamma$ be the quotient map. For representatives $\mathcal H_j$ of the $\Gamma$-orbits of branch mirrors, we define
\begin{equation}\label{eq:branch-definition}
 B_j=\overline{p_\Gamma(\mathcal H_j)}^{\,Y_\Gamma},\qquad
 B_\Gamma=\sum_jB_j.
\end{equation}
Write $T(\Gamma)$ for the normalized degree of this reduced divisor and $T(L)$ for the full case of the lattice group, that is, 
\begin{equation}\label{eq:Tintro}
 T(\Gamma)=\frac{B_\Gamma\cdot\lambda_\Gamma^{n-1}}{\lambda_\Gamma^n},
 \qquad T(L)=T(\PO^+(L)).
\end{equation}
For a mirror $\mathcal H$, its setwise stabilizer and pointwise stabilizer are
\[
 S_\Gamma(\mathcal H)=\{g\in\Gamma:g\mathcal H=\mathcal H\},\qquad
 I_\Gamma(\mathcal H)=\{g\in\Gamma:g|_{\mathcal H}=\mathrm{id}\}.
\]

\begin{lemma}\label{lem:inertia}
For a branch mirror $\mathcal H=\D_v$ with projective reflection $r=r_v$, we have
\[
 I_\Gamma(\mathcal H)=\langle r\rangle=\{1,r\},\qquad
 S_\Gamma(\mathcal H)=C_\Gamma(r).
\]
The generic point stabilizer on $\mathcal H$ is $\{1,r\}$, the ramification index is two, and $B_\Gamma$ has finitely many irreducible components.
\end{lemma}

\begin{proof}
A nonidentity projective orthogonal transformation fixing a divisor pointwise is its unique projective reflection. Hence
\[
 I_\Gamma(\mathcal H)=\{1,r\},\qquad
 g\mathcal H=\mathcal H\ \Longleftrightarrow\ grg^{-1}=r.
\]
At a generic point, the local quotient and its ramification index are
\[
 (z,w)\longmapsto(u,w),\quad u=z^2,\quad \mathcal H=\{z=0\}, \quad 
B:=\overline{p_\Gamma(\mathcal H)}^{\,Y_\Gamma}=\{u=0\},
\]
\begin{equation}\label{eq:ramification-index}
 e_\Gamma(\mathcal H)=\operatorname{ord}_{\mathcal H}(p_\Gamma^*u)=2.
\end{equation}
In particular, for nonzero local meromorphic functions or sections of a line bundle on $X_\Gamma$, we have
\[
 \operatorname{ord}_{\mathcal H}(p_\Gamma^*f)=2\operatorname{ord}_{B}(f).
\]

Choose a torsion-free normal subgroup $\Gamma'\triangleleft\Gamma$ of finite index, and let $\pi:Y_{\Gamma'}\to Y_\Gamma$ be the induced finite algebraic morphism. On the open modular varieties,
\[
 \D_L\xrightarrow{p_{\Gamma'}}X_{\Gamma'}\xrightarrow{\pi}X_\Gamma,
 \qquad p_\Gamma=\pi\circ p_{\Gamma'}.
\]
Since $\Gamma'$ acts freely on $\D_L$, the quotient map $p_{\Gamma'}$ is locally biholomorphic. On a sufficiently small neighborhood $U\subset\D_L$, the identity $\pi=p_\Gamma\circ(p_{\Gamma'}|_U)^{-1}$ therefore identifies the local ramification of $\pi$ with that of $p_\Gamma$. In particular, the local model above shows that each $p_\Gamma(\mathcal H_j)$ is a branch divisor of $\pi|_{X_{\Gamma'}}$. Conversely, since $\pi$ is finite, a branch divisor on $X_\Gamma$ comes from a ramification divisor on $X_{\Gamma'}$. Lifting locally through $p_{\Gamma'}$ gives a divisorial ramification component of $p_\Gamma$. By \cite[Proposition 1.2]{GH14}, in effective projective notation, every such component is a mirror $\D_v$ with $r_v\in\Gamma$. Thus the divisorial branch components on $X_\Gamma$ are exactly the images of branch mirrors.

The Baily--Borel boundary has codimension at least two, so every codimension-one branch component of $\pi$ meets $X_\Gamma$. Taking closures gives
\[
 \begin{gathered}
 \{\text{codimension-one irreducible components of }\operatorname{Br}(\pi)\}
 =\{B_j\}_j.
 \end{gathered}
\]
The non-\'etale locus of $\pi$ is Zariski closed, and a finite morphism sends closed subsets to closed subsets. Hence its image $\operatorname{Br}(\pi)$ is a closed algebraic subset of $Y_\Gamma$ \cite[Tags 02GH and 01WM]{Stacks}. Since $Y_\Gamma$ is Noetherian, this branch locus has finitely many irreducible components \cite[Tag 0052]{Stacks}. Its codimension-one components therefore form a finite list $B_1,\ldots,B_s$. Each is a prime algebraic divisor. Thus
\[
 B_\Gamma=\sum_{j=1}^sB_j,\qquad
 \operatorname{Supp}(B_\Gamma)=\bigcup_{j=1}^sB_j.
\]
The support is a closed algebraic subset of $Y_\Gamma$, being a finite union of such subsets.
\end{proof}

For a reflective form $F$, choose representatives $\mathcal H_1,\ldots,\mathcal H_t$ of its zero-mirror orbits and put
\begin{equation}\label{eq:zero-divisor-definition}
 B_j=\overline{p_\Gamma(\mathcal H_j)}^{\,Y_\Gamma},\qquad
 m_j=\operatorname{ord}_{\mathcal H_j}(F),\qquad
 B_F=\sum_{j=1}^tm_jB_j.
\end{equation}
The coefficients of $B_F$ are measured on the domain, before taking the quotient. If $m_j\leq N$ for all $j$, where $N\geq1$ is a fixed integer, then
\begin{equation}\label{eq:bounded-divisor}
 0\leq B_F\leq NB_\Gamma.
\end{equation}
The next lemma is the Baily--Borel form of the reflective-divisor calculation in \cite[Section 2, proof of Theorem 2.1]{GH14}. Its lower bound is the input needed in Proposition~\ref{prop:fixed-carrier}.

\begin{lemma}\label{lem:divisor}
Let $F$ be a reflective modular form of weight $k\in\frac12\Z_{>0}$ with a finite-order multiplier system for $\Gamma\subset\mathrm{O}(n,2)$ as above. Writing $\Cl(Y_\Gamma)$ for the Weil divisor class group, we have
\begin{equation}\label{eq:weighted-divisor}
 [B_F]=\sum_jm_j[B_j]=2k\lambda_\Gamma
 \quad\text{in }\Cl(Y_\Gamma)\otimes\Q.
\end{equation}
If every zero multiplicity is at most $N\geq1$, then
\begin{equation}\label{eq:divisor}
 n-2\leq2k\leq N\cdot T(\Gamma).
\end{equation}
The equality $2k=N\cdot T(\Gamma)$ holds if and only if $B_F=NB_\Gamma$.
\end{lemma}

\begin{proof}
Choose $q$ sufficiently divisible to kill the multiplier and make $qk$ a multiple of $a$. Then $F^q\in H^0(Y_\Gamma,\mathcal L_{\Gamma,qk})$, and Lemma~\ref{lem:inertia} gives $\operatorname{ord}_{B_j}(F^q)=qm_j/2$.
Reflectivity accounts for all interior zero divisors, and the Baily--Borel boundary contains no prime divisor. Therefore,
\begin{equation}\label{eq:power-divisor}
 \divisor_{Y_\Gamma}(F^q)=\frac q2B_F,\qquad
 qk\lambda_\Gamma=\frac q2[B_F].
\end{equation}
Divide by $q/2$ to obtain the class identity. Intersecting \eqref{eq:bounded-divisor} with $\lambda_\Gamma^{n-1}$ gives
\[
 2k\lambda_\Gamma^n
 =B_F\cdot\lambda_\Gamma^{n-1}
 \leq NB_\Gamma\cdot\lambda_\Gamma^{n-1}.
\]
Dividing by $\lambda_\Gamma^n>0$ leads to the inequality. The ampleness gives the equality criterion, and the singular-weight bound yields $2k\geq n-2$.
\end{proof}

For a reflective form, the intersection of $[B_F]=2k\lambda_\Gamma$ with $\lambda_\Gamma^{n-1}$ recovers the weight from the weighted zero divisor as
\begin{equation}\label{eq:weight-from-divisor}
 k=\frac{B_F\cdot\lambda_\Gamma^{n-1}}{2\lambda_\Gamma^n}.
\end{equation}
For reflective forms, this identity is the intersection-theoretic
counterpart of Bruinier's analytic divisor-degree formula
\cite[Theorem~1]{Bru04}. For a meromorphic modular form $F$ of
integral weight $k$ on $X_\Gamma$ by a freely
acting arithmetic group in the identity component, his formula reads
\[
\int_{\operatorname{div}(F)} \Omega^{n-1}
= k \int_{X_\Gamma} \Omega^n,
\]
where $\Omega$ is the first Chern form of the weight-one automorphic
line bundle. The formula applies for $n\geq 4$ and for $n=3$ when the rational Witt index is at most one. In our notation, the divisor on the quotient is $\tfrac12 B_F$. Lemma~\ref{lem:divisor} gives the
algebraic identity also for lattices of rational Witt index two
in dimension three. 

The numerical projection formula \cite[Lemma 33.45.7, Tag 0BET]{Stacks} says that, for a proper dominant generically finite map $f:Z\to W$ of $d$-dimensional integral projective varieties and a line bundle $\mathcal A$ on $W$, one has
\begin{equation}\label{eq:projection-formula}
 \bigl((f^*\mathcal A)^d\cdot Z\bigr)
 =\deg(f)\,(\mathcal A^d\cdot W).
\end{equation}

Fix an effective rational orthogonal arithmetic group $\Lambda<G$. 
For any subgroups $\Gamma\leq\Lambda$ of finite index, the compactified quotient map satisfies
\[
 p:Y_\Gamma\longrightarrow Y_\Lambda,\qquad
 \deg(p)=[\Lambda:\Gamma],\qquad p^*\lambda_\Lambda=\lambda_\Gamma.
\]
By \eqref{eq:projection-formula} in dimension $n$, we have $\lambda_\Gamma^n=[\Lambda:\Gamma]\lambda_\Lambda^n$.
Write
\[
 B_\Lambda=B_1+\cdots+B_s,\qquad
 c_j=\frac{B_j\cdot\lambda_\Lambda^{n-1}}{\lambda_\Lambda^n}>0.
\]
For each $j$, choose a mirror $\mathcal H_j$ above $B_j$ and its projective reflection $r_j\in\Lambda$. We use the fixed-point ratios for $\Lambda/\Gamma$ introduced above.

\begin{proposition}\label{prop:branch-cover}
For every finite-index subgroup $\Gamma\leq\Lambda$, we have
\begin{equation}\label{eq:branch-cover}
 T(\Gamma)=\sum_{j=1}^s c_j\cdot \fpr_\Gamma(r_j).
\end{equation}
In particular, $T(\Gamma)\leq T(\Lambda)$.
\end{proposition}

\begin{proof}
Let $p:Y_\Gamma\to Y_\Lambda$ be the finite cover of degree $[\Lambda:\Gamma]$. Choose a generic point $y$ of $B_j$ and a lift $z\in\mathcal H_j$ with stabilizer $\langle r_j\rangle$ in $\Lambda$, as in Lemma~\ref{lem:inertia}. The map 
\[
\Lambda/\Gamma \to p^{-1}(y), \quad g\Gamma\mapsto[g^{-1}z]_\Gamma
\]
is surjective and identifies $p^{-1}(y)$ with $\langle r_j\rangle\backslash(\Lambda/\Gamma)$, since $[g_1^{-1}z]_\Gamma=[g_2^{-1}z]_\Gamma$ iff $g_2\gamma g_1^{-1}\in \langle r_j\rangle$ for some $\gamma \in \Gamma$. 
The stabilizer of $g^{-1}z$ in $\Gamma$ is $\Gamma\cap g^{-1}\langle r_j\rangle g$. It is nontrivial if and only if $r_jg\Gamma=g\Gamma$, in which case $g^{-1}\HH_j$ is a branch mirror for $\Gamma$. When the stabilizer is trivial, $p_\Gamma$ is locally biholomorphic and the point of $[g^{-1}z]_\Gamma$ is not in $B_\Gamma$. It follows that
\[
|p^{-1}(y)\cap\operatorname{Supp}(B_\Gamma)|
=
|\Fix_{\Lambda/\Gamma}(r_j)|.
\]
For each branch component above $B_j$, the number of its points over a generic $y$ equals its degree over $B_j$. Every component of $B_\Gamma$ maps to some $B_j$, so the pushforward is
\begin{equation}\label{eq:branch-pushforward}
 p_*B_\Gamma
 =\sum_{j=1}^s|\Fix_{\Lambda/\Gamma}(r_j)|\,B_j.
\end{equation}
Using $p^*\lambda_\Lambda=\lambda_\Gamma$, $\lambda_\Gamma^n=[\Lambda:\Gamma]\lambda_\Lambda^n$, and \eqref{eq:projection-formula} on each branch component gives
\[
T(\Gamma)=\frac{B_\Gamma\cdot(p^*\lambda_\Lambda)^{n-1}}
         {[\Lambda:\Gamma]\lambda_\Lambda^n}
=\frac{(p_*B_\Gamma)\cdot\lambda_\Lambda^{n-1}}
         {[\Lambda:\Gamma]\lambda_\Lambda^n}
=\sum_{j=1}^sc_j\cdot\fpr_\Gamma(r_j)\leq T(\Lambda). \qedhere
\]
\end{proof}

\begin{corollary}\label{cor:fixed-carrier-decay}
If the finite-index subgroups $\Gamma_i\leq\Lambda$ satisfy $[\Lambda:\Gamma_i]\to\infty$, then $T(\Gamma_i)\to0$.
\end{corollary}

\begin{proof}
The sum in \eqref{eq:branch-cover} has finitely many terms. Apply Proposition~\ref{prop:fpr} to each $r_j\neq1$.
\end{proof}

\begin{proposition}\label{prop:fixed-carrier}
Fix an integer $N\geq1$ and an integral lattice $L$ of signature $(n,2)$, with $n\geq3$. Then modulo $\C^\times$ there are only finitely many reflective modular forms for all finite-index subgroups of $\PO^+(L)$ whose zero multiplicities are at most $N$. 
\end{proposition}

\begin{proof}
If the subgroup indices were unbounded, then we could choose reflective modular forms of weights $k_i$ with zero multiplicities at most $N$ for subgroups $\Gamma_i\leq\PO^+(L)$ with $[\PO^+(L):\Gamma_i]\to\infty$. Corollary~\ref{cor:fixed-carrier-decay} and Lemma~\ref{lem:divisor} would give
\[
n-2\leq2k_i\leq N\cdot T(\Gamma_i)\longrightarrow0,
\] 
a contradiction. By \cite[Theorem 6.12]{BH62}, $\PO^+(L)$ is finitely generated, so it has only finitely many subgroups of each bounded index. This also follows from \cite[Theorem 8.1]{Wan83}. 

Fix one of such groups $\Gamma$. By Lemma~\ref{lem:inertia}, $B_\Gamma$ has finitely many irreducible components, say $s$. Each coefficient of $B_F$ lies in $\{0,\ldots,N\}$, so there are at most $(N+1)^s$ choices for $B_F$. If $F_1$ and $F_2$ have the same $B_F$, then $F_1/F_2$ is constant by the Koecher principle. 
\end{proof}

\section{Effective bounds for branch degrees}\label{sec:arithmetic}
In this section, we adapt Ma's approach \cite{Ma18} to establish an upper bound for $T(L)=T(\PO^+(L))$. For an integral lattice $M$ of rank $n+2$, fix its dual lattice, discriminant group, and discriminant exponent as follows
\[
 M^\vee=\{x\in M\otimes\Q:\, (x,M)\subseteq\Z\},\qquad
 A_M=M^\vee/M,\qquad D(M)=\exp(A_M).
\]
The exponent is the least positive integer annihilating $A_M$. If $d_1\mid\cdots\mid d_r$ are the Smith invariant factors of a Gram matrix of $M$, then
\begin{equation}\label{eq:smith}
 \cont(M):=\gcd\{(x,y):x,y\in M\}=d_1,\quad D(M)=d_r,\quad |\det M|=\prod_{i=1}^rd_i.
\end{equation}
The lattice $M$ is called primitive if $\cont(M)=1$, in which case
\begin{equation}\label{eq:detbound}
 |\det M|\leq D(M)^{n+1}.
\end{equation}
Write $M(a)$ for $M$ with its bilinear form multiplied by $a>0$. This scaling preserves the orthogonal group, the symmetric domain, and the reflections. In particular, dividing the bilinear form of $M$ by $\cont(M)$ gives a primitive integral lattice with the same orthogonal group and mirrors.

We will estimate $T(L)$ for a primitive integral lattice $L$ of signature $(n,2)$ with $n\geq 3$. Let $\Delta$ be a rational orthogonal arithmetic group on a type-IV domain $\mathcal D$ of dimension $d$ with $d\geq 2$, and use $Y_\Delta$ and $\lambda_\Delta$ as in Section~\ref{sec:main step}. With compatible invariant volume forms on $\mathcal D$ and its compact dual $\mathcal D^{(c)}$, the Hirzebruch--Mumford volume is (see \cite[Section 2, equation (5)]{GHS07})
\[
 v(\Delta):=\vol_{\HM}(\Delta)
 =\frac{\vol(\Delta\backslash\mathcal D)}{\vol(\mathcal D^{(c)})},
\]
where $\mathcal D^{(c)}$ is the smooth projective quadric obtained by dropping the inequality in the domain definition. The quotient volume is for the effective action. In each fixed dimension, $\vol_{\HM}$ is the ordinary invariant quotient volume times a positive constant; ambient domains and mirrors use the constants in their respective dimensions.

Along sufficiently divisible even weights, \cite[Proposition 1.2 and its proof]{GHS07}, equivalently \cite[Section 4.2, equation (4.1)]{Ma18}, gives
\begin{equation}\label{eq:HM-asymptotic}
 \dim M_m(\Delta)=\frac{2v(\Delta)}{d!}m^d+O(m^{d-1}).
\end{equation}
Here, $M_m(\Delta)$ has trivial character and, when $d=2$, includes holomorphy at the cusps. For a sufficiently divisible even weight $a$, the Hilbert polynomial of the ample bundle $\mathcal L_{\Delta,a}$ gives 
\[
 \dim M_{aj}(\Delta)
 =h^0(Y_\Delta,\mathcal L_{\Delta,a}^{\otimes j})
 =\frac{\lambda_\Delta^d}{d!}(aj)^d+O(j^{d-1}).
\]
Comparison of leading coefficients yields
\begin{equation}\label{eq:HM-intersection}
 \lambda_\Delta^d=2v(\Delta).
\end{equation}
For a finite-index inclusion $\Delta'\subset\Delta$ before projectivization, \cite[Section 4.2, equation (4.2)]{Ma18} gives
\begin{equation}\label{eq:HM-index}
 v(\Delta')=[\langle\Delta,-1\rangle:\langle\Delta',-1\rangle]v(\Delta).
\end{equation}
The index is equal to $[\Delta:\Delta']$ if both groups contain $-1$, or if both are already effective.

\begin{lemma}\label{lem:mumford-degrees}
Let $\Gamma$ be an effective rational orthogonal arithmetic group on $\D_L$. For a branch mirror $\mathcal H$, let $B$ be the reduced closure of its image in $Y_\Gamma$. Using the stabilizers in Lemma \ref{lem:inertia}, its effective mirror group is 
\[
 \Delta_{\mathcal H}:=\overline S_\Gamma(\mathcal H)
 =S_\Gamma(\mathcal H)/I_\Gamma(\mathcal H).
\]
Then
\begin{equation}\label{eq:ambient-mirror-degrees}
 \lambda_\Gamma^n=2v(\Gamma),\qquad
 B\cdot\lambda_\Gamma^{n-1}
 =\lambda_{\Delta_{\mathcal H}}^{n-1}=2v(\Delta_{\mathcal H}).
\end{equation}
\end{lemma}

\begin{proof}
The ambient identity is \eqref{eq:HM-intersection}. The mirror quotient normalizes the open branch component by \cite[Section 4.1, before Lemma 4.2]{Ma18} and \cite[Appendix A, Propositions A.5--A.6]{YZ20}. On Baily--Borel compactifications, it gives a finite birational map
\[
 \nu:Y_{\Delta_{\mathcal H}}\longrightarrow B,
 \qquad\deg(\nu)=1.
\]
The mirror inclusion preserves the automorphic bundle, including its compactified extension. For a common sufficiently divisible even weight $a$,
\begin{equation}\label{eq:mirror-line-pullback}
 \nu^*(\mathcal L_{\Gamma,a}|_B)
 \simeq\mathcal L_{\Delta_{\mathcal H},a},\qquad
 \nu^*(\lambda_\Gamma|_B)=\lambda_{\Delta_{\mathcal H}};
\end{equation}
see \cite[Proposition 6.25 (b)]{Pin90} and \cite[proof of Theorem A.14]{YZ20} for functoriality of the compactified bundles. Applying \eqref{eq:projection-formula} to $\nu$ with $\mathcal A=\mathcal L_{\Gamma,a}|_B$ and dividing by $a^{n-1}$, we derive 
\[
 \begin{aligned}
 \deg(\nu)\,(B\cdot\lambda_\Gamma^{n-1})
 &=\frac{\bigl((\nu^*(\mathcal L_{\Gamma,a}|_B))^{n-1}
                    \cdot Y_{\Delta_{\mathcal H}}\bigr)}{a^{n-1}}\\
 &=\lambda_{\Delta_{\mathcal H}}^{n-1}
 =2v(\Delta_{\mathcal H}).
 \end{aligned}
\]
Since $\deg(\nu)=1$, this proves the mirror identity, also when the mirror has dimension two.
\end{proof}

For a branch mirror, Lemma~\ref{lem:inertia} gives $I_\Gamma(\mathcal H)=\langle r\rangle$. This kernel is removed in $\Delta_{\mathcal H}$, so it contributes no further factor two. Substituting \eqref{eq:ambient-mirror-degrees} into the definition of $T$ gives
\begin{equation}\label{eq:T-mirror-HM}
 T(\Gamma)=\sum_{B\subset B_\Gamma}
       \frac{v(\overline S_\Gamma(\mathcal H_B))}{v(\Gamma)},
\end{equation}
where $\mathcal H_B$ is a branch mirror above $B$. 

Let $l\in L$ have positive norm and $\sigma_l\in\OO^+(L)$. Take $l$ primitive, so $\Z l=L\cap\Q l$. Put
\[
 K_l=l^\perp\cap L,\qquad
 \Gamma_l=\im\bigl(\{g\in\OO^+(L):g(l)=l\}\longrightarrow\OO^+(K_l)\bigr).
\]
Its effective image is the mirror group of Lemma~\ref{lem:mumford-degrees}, and
\[
 -1=(-\sigma_l)|_{K_l}\in\Gamma_l,\qquad
 \Gamma_l/\{\pm1\}\simeq
 \overline S_{\PO^+(L)}(\D_l).
\]
Set $B_L=B_{\PO^+(L)}$ and $\lambda_L=\lambda_{\PO^+(L)}$. For the component $B_l$ associated with $l$, Lemma~\ref{lem:mumford-degrees} gives
\[
 \lambda_L^n=2v(\OO^+(L)),\qquad
 B_l\cdot\lambda_L^{n-1}=\lambda_{\Gamma_l}^{n-1}=2v(\Gamma_l).
\]
Since $T(L):=T(\PO^+(L))$ and $B_L=\sum_{[l]}B_l$, \eqref{eq:T-mirror-HM} becomes
\begin{equation}\label{eq:T-HM}
 T(L)=\sum_{[l]}\frac{v(\Gamma_l)}{v(\OO^+(L))}.
\end{equation}
The sum is over $\OO^+(L)$-orbits of primitive reflective vectors. The reflection $\sigma_l$ puts $l$ and $-l$ in the same orbit. For an integral lattice $M$, define its split sum using $K_l=l^\perp\cap M$,
\begin{equation}\label{eq:split-sum}
 S(M)=\sum_{\substack{[l]\bmod\OO^+(M)\\l\text{ primitive},\ (l,l)>0\\M=\Z l\perp K_l}}
             \frac{v(\OO^+(K_l))}{v(\OO^+(M))}.
\end{equation}
We use two results of Ma to reduce to split reflective vectors. Both apply for $n\geq3$, with the bilinear form reversed for Ma's signature convention $(2,n)$. The criterion in \cite[Lemmas 4.1--4.2]{Ma18} distinguishes split and non-split reflective vectors.

\begin{lemma}\label{lem:Ma-reflective}
Let $L$ have signature $(n,2)$, with $n\geq3$. A primitive positive-norm vector $l$ is reflective, that is, $\sigma_l\in \mathrm{O}^+(L)$ if and only if
\[
 [L:\Z l\perp K_l]\in\{1,2\}.
\]
In the split case, this index is one and $\Gamma_l=\OO^+(K_l)$. In the non-split case, $\Gamma_l$ is the stabilizer of an element of order two in $A_{K_l}$ and
\[
 [\OO^+(K_l):\Gamma_l]\leq2^{n+1}.
\]
\end{lemma}

Let $R_{\mathrm{ns}}(L)$ denote the non-split reflective orbits, and put
\[
 \mathcal N(L)=\sum_{[l]\in R_{\mathrm{ns}}(L)}
          \frac{v(\OO^+(K_l))}{v(\OO^+(L))}.
\]
Let $\mathcal T$ be the set of actual sublattices $M_l=\Z l\perp K_l$ arising from these vectors. For $M\in\mathcal T$, its stabilizer is
\[
 H_M=\{g\in\OO^+(L):gM=M\}=\OO^+(L)\cap\OO^+(M).
\]
Ma's orbit estimate \cite[Lemma 6.5]{Ma18} bounds the non-split sum in terms of split sums. 

\begin{lemma}\label{lem:Ma-orbits}
For an integral lattice $L$ of signature $(n,2)$ with $n\geq3$, one has
\[
\mathcal N(L)\leq\sum_{[M]\in\mathcal T/\OO^+(L)} [\OO^+(L):H_M]\cdot S(M).
\]
\end{lemma}

The index counts the sublattices in the orbit of $M$, and $S(M)$ is constant on that orbit. Hence
\[
 \mathcal N(L)\leq\sum_{M\in\mathcal T}S(M)
 \leq\sum_{\substack{M\subset L\\{}[L:M]=2}}S(M).
\]
By Lemma~\ref{lem:Ma-reflective} and \eqref{eq:HM-index},
\[
 \begin{aligned}
 T(L)
 &=S(L)+\sum_{[l]\in R_{\mathrm{ns}}(L)}
       [\OO^+(K_l):\Gamma_l]
       \frac{v(\OO^+(K_l))}{v(\OO^+(L))}\\
 &\leq S(L)+2^{n+1}\mathcal N(L).
 \end{aligned}
\]
Consequently,
\begin{equation}\label{eq:reduction}
 T(L)\leq S(L)+2^{n+1}
             \sum_{\substack{M\subset L\\{}[L:M]=2}}S(M).
\end{equation}
The sum has $2^{n+2}-1$ terms, corresponding to the kernels of nonzero maps $L/2L\to\F_2$.

For a prime $p$, let $\Z_p$ be the ring of $p$-adic integers and $v_p$ the exponent of $p$ in a positive integer. Use a Jordan decomposition over $\Z_p$ with $M_{p,j}$ either zero or unimodular, meaning that its Gram determinant is a unit in $\Z_p$. With the scaling convention above,
\[
 M\otimes\Z_p=\mathop{\perp}_{j\geq0}M_{p,j}(p^j),\qquad
 r_{p,j}=\rk M_{p,j},\qquad
 \mu_p=\max\{j:r_{p,j}>0\}=v_p(D(M)).
\]
Set
\[
 m_{p,j}(M)=\sum_k|k-j|r_{p,k}-\mu_p,\qquad
 \mathcal P(M)=\{p>2:\, p\mid D(M),\ r_{p,j}\leq1\text{ for all }j\}.
\]
Following \cite[Definitions 5.4 and 5.7]{Ma18}, for $r_{p,j}>0$, define
\[
 \eps_{p,j}(M)=
 \begin{cases}
 2^{-m_{2,j}(M)/2},&p=2,\\
 4p^{-m_{p,j}(M)/2},&p\in\mathcal P(M),\\
 (1+p^{-\lfloor r_{p,j}/2\rfloor})p^{-m_{p,j}(M)/2},
       &p>2,\ p\notin\mathcal P(M),
 \end{cases}
\]
where $\lfloor x\rfloor:=\max\{ m\in\Z:\, m\leq x \}$, and set
\[
 \eps_p(M)=\sum_{j:\, r_{p,j}>0}\eps_{p,j}(M),\qquad
 \eps(M)=\prod_{p\mid2D(M)}\eps_p(M).
\]
Since $\eps_2(M)=1$ when $D(M)$ is odd, this agrees with Ma's product over $p\mid D(M)$.

For primitive lattices, the local exponents satisfy Ma's bound \cite[Lemma 5.9]{Ma18}.

\begin{lemma}\label{lem:Ma-local-exponents}
For a primitive lattice $M$ of rank $n+2$ and every nonzero Jordan block,
\begin{equation}\label{eq:local-exponents}
 m_{p,j}(M)\geq\max\{0,\, n-r_{p,j}\}.
\end{equation}
\end{lemma}

\begin{proof}
For every $n\geq3$, primitivity gives $r_{p,0}>0$. The claim is immediate if $\mu_p=0$. Otherwise, subtract one from each endpoint rank, writing $\delta_{i,j}=1$ when $i=j$ and $0$ otherwise. From
\[
\widehat r_k:=r_{p,k}-\delta_{k,0}-\delta_{k,\mu_p} \quad \text{and} \quad
|0-j|+|\mu_p-j|=\mu_p,
\]
we conclude that for $0\leq j\leq\mu_p$,
\begin{align*}
 m_{p,j}(M)&=\sum_k|k-j|\widehat r_k\geq0,\\
 m_{p,j}(M)&\geq\sum_{k\ne j}\widehat r_k
 =n-r_{p,j}+\delta_{j,0}+\delta_{j,\mu_p}\geq n-r_{p,j}.\qedhere
\end{align*}
\end{proof}

The split-sum estimate that we use is \cite[Lemma 6.3 and equation (6.2)]{Ma18}. 

\begin{lemma}\label{lem:Ma-split}
For every integral $M$ of signature $(n,2)$ with $n\geq3$, one has
\begin{equation}\label{eq:Ma-split}
 S(M)\leq f_n\eps(M)D(M)^{-1/2},\qquad
 f_n=\frac{2^7\cdot9\,\pi^{n/2+1}\zeta(\lfloor n/2\rfloor+1)}
                {\Gamma(n/2+1)}.
\end{equation}
Here, $\zeta$ and $\Gamma(\,\cdot\,)$ are the Riemann zeta and Euler gamma functions. The estimate applies to primitive and nonprimitive $M$.
\end{lemma}

We first estimate the split sum on a primitive lattice. For an index-two sublattice we use the same estimate after dividing its bilinear form by its content.

\begin{lemma}\label{lem:coarse-local}
Let $L$ be primitive of rank $n+2$. Set $D=D(L)$ and $W_n(D)=(4n+8)^{\omega(2D)}$, where $\omega(t)$ counts the distinct prime divisors of $t$. Then
\[
 S(L)\leq f_nW_n(D)D^{-1/2},\qquad
 S(M)\leq 2f_nW_n(D)D^{-1/2}\quad([L:M]=2).
\]
\end{lemma}

\begin{proof}
Lemma~\ref{lem:Ma-local-exponents} yields $m_{p,j}(L)\geq0$. Each nonzero Jordan block contributes at least one to $\sum_jr_{p,j}=n+2$. It follows that
\[
 \eps_2(L)\leq n+2,\qquad \eps_p(L)\leq 4n+8 \quad (p>2),\qquad
 \eps(L)\leq W_n(D).
\]
The estimate for $S(L)$ follows from \eqref{eq:Ma-split}.

If $[L:M]=2$, then $2L\subset M$, so $c=\cont(M)$ divides $\cont(2L)=4$. Put $M_0=M(1/c)$, which is primitive. Scaling preserves the split vectors, their orbits, and the volume ratios in \eqref{eq:split-sum}, so $S(M)=S(M_0)$. Since $L^\vee\subset M^\vee$ and $D(M)M^\vee\subset M\subset L$, we have $D(L)\mid D(M)$. Scaling the Smith invariant factors gives $D(M_0)=D(M)/c\geq D/4$. At odd primes, $M$ and $L$ have the same local lattice and $c$ is a unit. Thus $D(M_0)$ and $D(L)$ have the same odd prime divisors, giving $W_n(D(M_0))=W_n(D)$. The primitive estimate gives
\[
 S(M)=S(M_0)\leq f_n W_n(D)(D/4)^{-1/2}=2f_n W_n(D)D^{-1/2}. \qedhere
\]
\end{proof}

Combining these estimates gives the effective branch-degree bound.

\begin{theorem}\label{thm:bounds}
For each $n\geq3$, there is an effective constant $C_n>0$ such that every primitive integral lattice $L$ of signature $(n,2)$ satisfies
\begin{equation}\label{eq:coarse-intro}
 T(L)\leq C_nD(L)^{-1/4}.
\end{equation}
One may take
\begin{equation}\label{eq:Cn}
 C_n=f_n(1+4^{n+2})b_n^{b_n^4},\qquad b_n=4(n+2).
\end{equation}
\end{theorem}
\begin{proof}
Splitting the prime divisors of $2D(L)$ into those at most $b_n^4$ and those larger, and using $b_n<p^{1/4}$ for the larger primes, we obtain
\[
 b_n^{\omega(2D(L))}\leq b_n^{b_n^4}D(L)^{1/4}.
\]
Substitute the bounds from Lemma~\ref{lem:coarse-local} into \eqref{eq:reduction}. The number of index-two sublattices gives
\begin{align*}
 T(L)&\leq f_nW_n(D(L))D(L)^{-1/2}
            [1+2^{n+2}(2^{n+2}-1)]\\
 &\leq f_n(1+4^{n+2})b_n^{\omega(2D(L))}D(L)^{-1/2}\\
 &\leq C_nD(L)^{-1/4}.\qedhere
\end{align*}
\end{proof}

Note that Ma \cite{Ma18} obtained the stronger bound $T(L)\ll_n D(L)^{-1/2}$ for any $n\geq 4$; our estimate $T(L)\ll_n D(L)^{-1/4}$ holds for all $n\geq 3$.

\section{Finiteness of reflective modular forms}\label{sec:finiteness}
In this section, we complete the proof of Theorem \ref{Main theorem}. 

\begin{proposition}\label{prop:finite-carriers}
Fix an integer $N\geq 1$. Up to isometry, only finitely many primitive integral lattices $L$ admit a reflective form for a finite-index subgroup of $\PO^+(L)$ whose zero multiplicities are at most $N$. 
\end{proposition}

\begin{proof}
Let $F$ be reflective for $\Gamma\leq\PO^+(L)$, with zero multiplicities at most $N$. For each zero mirror, $r_v\in\Gamma$ implies $\sigma_v\in\OO^+(L)$ since $-1\in\OO^+(L)$. A sufficiently divisible power of $F$ is a reflective form of integral weight for the inverse image of $\Gamma$ in $\OO^+(L)$. If $L$ is odd, use the even lattice $L(2)$, which has the same orthogonal group and reflections. Assume that the signature of $L$ is $(n,2)$. \cite[Theorem 1.2]{Wan25b} yields $n\leq 20$ or $n=22, 26$. 

For the original primitive $L$, Lemma~\ref{lem:divisor}, Proposition~\ref{prop:branch-cover}, and Theorem~\ref{thm:bounds} give
\begin{equation}\label{eq:main-inequalities}
1\leq n-2\leq2k\leq N\cdot T(\Gamma)\leq N\cdot T(L)
 \leq NC_nD(L)^{-1/4}.
\end{equation}
It follows that $D(L)\leq (NC_n)^4$. By \eqref{eq:detbound}, there are only finitely many determinant values, and \cite[Theorem 6.10]{SP23} gives finitely many isometry classes in each rank. 
\end{proof}

\begin{proof}[Proof of Theorem~\ref{Main theorem}]
Fix an integer $N\geq1$. Let $F$ be a reflective modular form for $\Gamma\leq\PO^+(L)$, with every zero multiplicity at most $N$. We set
\[
 c=\cont(L),\qquad L_0=L(1/c),\qquad
 (x,y)_{L_0}=c^{-1}(x,y)_L.
\]
Then $L_0$ is integral and primitive. Passing from $L$ to $L_0$ leaves the symmetric domain, orthogonal group, and reflections unchanged, so $F$ is reflective for the same $\Gamma$ on $L_0$, with the same zero multiplicities. By Proposition~\ref{prop:finite-carriers}, there are finitely many such $L_0$ up to isometry. For lattices in a given isometry class, reflective modular forms coincide up to conjugation. For a fixed lattice $L_0$, Proposition~\ref{prop:fixed-carrier} shows that there are  finitely many such $F$ up to a nonzero scalar. This proves the finiteness of reflective modular forms with bounded zero multiplicities in full generality. 
\end{proof}

\smallskip
\noindent
\textbf{Acknowledgments.} 
H. Wang thanks Richard Borcherds, Valery Gritsenko, Shouhei Ma, Nils Scheithauer, Kaiwen Sun, and Brandon Williams for many helpful discussions related to this project over the past few years.

\smallskip

\noindent
\textbf{Disclosure of AI use.} 
The authors developed the initial mathematical framework and proof ideas. They subsequently used generative AI interactively in mathematical exploration, including literature searches and the development, testing, and refinement of parts of the proof. The authors independently verified the final arguments and take full responsibility for the content of the paper.

\bibliographystyle{plainnat}
\bibliofont
\bibliography{refs}

\end{document}